\documentclass[12pt]{article}
\usepackage{amsmath,epsfig,epsf,psfrag}
\newcommand{\LL}{{\cal L}}
\newcommand{\Ed}{{\rm Ed}}
\newcommand{\Fp}{{\rm Pra}}
\newcommand{\bgamma}{\mbox{\boldmath $\gamma$}}
\newcommand{\bbeta}{\mbox{\boldmath $\beta$}}
\newcommand{\boeta}{\mbox{\boldmath $\eta$}}
\newcommand{\balpha}{\mbox{\boldmath $\alpha$}}
\newcommand{\bA}{\mbox{\boldmath $A$}}
\begin{document}


\setlength{\textheight}{575pt} \setlength{\baselineskip}{23pt}

\title {A Latent Process Model for Dementia and Psychometric Tests}
\author{Julien Ganiayre$^{1}$, Daniel Commenges $^{1,2}$ and Luc Letenneur $^{2,3}$}
\date{}

\maketitle
1 INSERM, U875, Bordeaux, F-33076, France

2 Univ Bordeaux 2, Bordeaux, F33076, France

3 INSERM, U593, Bordeaux, F33076, France

\begin{abstract}
We jointly model longitudinal values of a psychometric test and diagnosis of dementia. The model is based on a continuous-time latent process representing cognitive ability. The link between the latent process and the observations is modeled in two phases. Intermediate variables are noisy observations of the latent process; scores of the psychometric test and diagnosis of dementia are obtained by categorizing these intermediate variables. We propose maximum likelihood inference for this model and we propose algorithms for performing this task. We estimated the parameters of such a model using the data of the five-year follow-up of the PAQUID study. In particularThis analysis yielded interesting results about the effect of educational level on both latent cognitive ability and specific performance in the mini mental test examination. The predictive ability of the model is illustrated by predicting diagnosis of dementia at the eight-year follow-up of the PAQUID study bsed on the information of the first five years.
\end{abstract}
\vspace{5mm} Key words: latent process, Brownian motion, joint
model, ordinal data, multivariate data, dementia, Alzheimer's disease, prediction.
\clearpage
\newpage

         \section{Introduction}

Alzheimer's disease is clinically characterized by a
progressive decline of cognitive abilities and is the main cause of dementia. This feature has two
important consequences for the modeling. First it is only an
idealization to consider that the disease starts at a particular
moment. The diagnosis is made at the time of examination by a
neurologist but this does not mean that the disease started at
this precise moment, nor even at any precise moment before
examination. The second consequence is that psychometric tests
which measure cognitive abilities can provide important information
regarding the progression of a pathological process which may lead to a
diagnosis of Alzheimer's disease or dementia. It is interesting to
devise models which link the two types of information (diagnosis
of dementia and psychometric tests) with three main objectives:
to better understand this link, to increase the power for detecting
risk factors, to predict dementia using previous observations of
scores of psychometric tests.

The problem can be tackled through joint modeling of an event (onset of dementia) and a
longitudinal marker (scores of a psychometric test). Joint modeling of CD4 cell counts and onset of AIDS or
death has been proposed by Faucett and Thomas (1996) and Wulfsohn and Tsiatis (1997). Concerning dementia
a model has been proposed by Jacqmin-Gadda, Commenges and Dartigues (2005), with the specific aim of estimating
a change-point in the regime of cognitive decline. Approaches based on a stochastic process framework are
particularly well suited to grasp the dynamics of diseases. Henderson, Diggle and Dobson (2000) proposed a model
in which a latent process acts as a time-dependent variable in a proportional hazards model. Other approaches of
joint modeling represent the event as the crossing of a barrier by the latent process (Whitmore, Crowder and
Lawless, 1998; Lee, DeGruttola and Schoenfeld, 2000). This approach was developed by Hashemi, Jacqmin-Gadda
and Commenges (2003) and applied to joint modeling of dementia and a psychometric test: in this model the
latent process was interpreted as representing cognitive ability. The present paper proposes an evolution of
this work with important differences which make the model much more flexible, and thus more usable; in
particular, for technical reasons, the Hashemi-Jacqmin-Gadda-Commenges model was restricted to linear time-trend
for the latent process.

We propose a new model which allows for the diagnosis of dementia and scores on a
psychometric test to be analyzed together. The
model looks particularly non-standard for dementia because we do not model onset of dementia but diagnosis of
dementia at the time of visit. This is in fact more realistic (even though interval-censoring was treated in the
Hashemi-Jacqmin-Gadda-Commenges model) because onset of dementia is an abstraction; cognitive decline is in fact 
most often progressive. Thus our basic model is that a neurologist makes a diagnosis of dementia if the subject
has a latent process below a certain threshold at the time of visit. As for scores of the psychometric test, we
consider a grid of threshold values $c_m$, and the subject has score $m$ if his latent process falls between
$c_m$ and $c_{m+1}$ at time of visit. This is a refined model compared with previous works treating ordinal
scores as continuous. With this approach, both diagnosis of dementia and score of the psychometric test are
categorized observations of the latent process. This is reminiscent of probit models for ordinal data (McCullagh
and Nelder, 1989; Chib and Greenberg, 1998), but here the underlying latent process allows to capture the dynamics
of the phenomenon under study. Our model is in fact slightly more complicated than the above description, as will be described later.

In section 2 we present a general form of the model which could be applied to contexts other than cerebral ageing. In section 3 the identifiability is studied and the
likelihood is derived. In section 4 we come to the
specific model used for dementia and the Mini Mental Score
Examination: we begin by describing the PAQUID study, a
large cohort study on ageing which provides the data we used; then we
describe the model, present a small simulation and give
results, particularly on the predictive ability of the model.
We end with a short conclusion.

         \section{Model and observations}

\subsection{Outline of the model}
We propose a general model for multidimensional longitudinal data based on a latent process. The observation of type $k$ for subject $i$ at time $t_{ij}$ will be denoted $Y_{ij}^k$ (in our application we will use observations of two types: $k=1$: diagnosis of dementia, $k=2$:  a psychometric test). Similarly as in Dunson (2003) we propose a hierarchical structure where the observations $Y_{ij}^k$ are possibly coarsening transformations of latent variables $\theta_{ij}^k$, and these latent variables are related to common latent elements.

The common latent element in our model is a latent process $\Lambda_i(t)$ which is defined in continuous time (in contrast with Dunson's model). In our
application it is natural to consider that subjects have a certain cognitive ability quantitatively represented
by $\Lambda_i(t)$ for any $t$, not only at measurement times. It is also made possible for this approach to treat unequally
spaced observation times which may be different from one subject to another. The model for the latent process,
driven by a Brownian motion, yields a natural correlation structure for the intermediate latent variables
$\theta_{ij}^k$, without introducing additional parameters which  would have to be estimated.

Another trait of our model is that it may be non-linear in the parameters. In the next section we present the model in the most general form that can be easily treated with our approach because it preserves the normality of the $\theta_{ij}^k$. Finally the model is a kind of multivariate probit model (Chib and Greenberg, 1998): it has a more direct interpretation than assuming that the $\theta_{ij}^k$ are related to the canonical parameters of a distribution in the exponential family, and it is related to threshold models already used by Hashemi, Jacqmin-Gadda and Commenges (2003) in this application. Moreover it leads to simpler numerical integrals.

Because of the central role of the latent process in our model, we will start by describing it, specifying
afterwards how it can be observed. We consider that there might be other observations (for instance other
psychometric tests) at other times; this would not affect our latent process which has an intrinsic meaning.

        \subsection{Latent process}

For each subject $i$ we introduce $ \Lambda_{i}=(\Lambda_{i}(t))_{t\ge 0}$, a continuous-time
stochastic latent process; in our application $ \Lambda_{i}(t)$ will represent the global cognitive ability of subject $i$ at time $t$. This latent process is modeled as a function of explanatory variables as:
\begin{equation} \label{latent}
 \Lambda_{i}(t) = f(\beta , x_{i}(t)) + F(\gamma,z_{i}(t)) a_{i} + W_i(t),
\end{equation}
where $W_i =(W_i(t))_{t\ge 0} $ is a standard Brownian motion. The $q$-vector of random effects $a_{i}$  has a
multivariate normal distribution: $a_{i} \sim \mathcal{N} (0,A)$; $a_{i}$ and $W_i$ are independent and the sets
$(a_i, i=1,\ldots,n)$ and  $(W_i , i=1,\ldots,n)$ are  sets of independent random vectors and processes;
the functions $f(.,.)$: $ {R}^{p} \times {R}^{l} \rightarrow {R}$ and $F(.,.)$: $ {R}^{p} \times {R}^{l}\rightarrow {R}^{q}$
are differentiable and possibly non-linear; $\beta$ and $\gamma$ are vectors of coefficients (some
which may be interpreted as regression coefficients, others which are used to parameterize the non-linearity) and
$x_{i}(t)$ and $z_i(t)$  are vectors of time dependent covariates including $t$ itself.

A linear model for the latent process $\Lambda_{i} (t)= x_{i}(t)^T\beta + z_{i}(t)^T a_{i} + W_i(t),$ is a particular case of model (\ref{latent}). Note that in a linear model there is no parameter $\gamma$.

In the application we might consider the  non-linear model:
$\Lambda_{i}(t) = \beta_1+\beta_2 x_{i2}+(\beta_3+\beta_4x_{i2}) x_{i1}(t)^{\beta_5} + a_{i1}+ W_i(t)$, where $x_{i1}(t)=t$ is time itself, $x_{i2}$ represents educational level. This model is non-linear in time, but also in the parameters; parametrizing the power of time ($\beta_5$)  offers more flexibility in modeling the effect of time.

        \subsection{Observation equations.}
We consider that the values of ``tests'' at different time points are indirect observations of the latent process; in our application the ``tests'' include both  psychometric tests and diagnosis of dementia. We model the link between the latent process and the tests in two phases:
first we introduce, for subject $i$, intermediate random variables $\theta_{ij}^k$ which can be seen as potential measurements for each test
$k=1,\ldots,K $ of $\Lambda_i(t_{ij})$; secondly we represent the values of the tests as functions of these intermediate variables. The reason for differentiating these two phases is that the $\theta_{ij}^k$ are linear in $\Lambda_{i}(t_{ij})$ and  have normal distributions while the tests functions may be non-linear and discontinuous. The times $t_{ij}$ will be treated as deterministic. They might be random but under the condition that the mechanism leading to incomplete data is ignorable, a condition under which  the likelihood treating these times as fixed leads to the same inference as the correct likelihood. We make the same assumption for possibly missing data.

    \subsubsection{ Definition of $\theta_{ij}^k$.}

The intermediate variables for subject $i$ and for test $k$ are defined as:
\begin{equation} \label{intermediate}
\theta_{ij}^k = \Lambda_{i}({t_{ij}}) +
g^{k}(\beta^k,x^k_{i} (t_{ij}) ) + G^{k}(\gamma^k,z^k_{i} (t_{ij}))
d_{i}^{k} + \epsilon^{k}_{ij},
\end{equation}
for $j=1,\ldots, n_i$,
where $g^{k}(.,.)$ and $G^{k}(.,.)$ are analogous to $f(.,.)$ and
$F(.,.)$ in the definition of the latent process but are specific to the $k^{th}$
test; $d_{i}^{k}$ is a $r_{k}$-random vector with normal distribution:
$d_{i}^{k} \sim \mathcal{N} (0,D^{k})$; the measurement errors
$\epsilon^{k}_{ij}$ are identically independently distributed (i.i.d) variables with normal distributions: $\epsilon^{k}_{ij} \sim \mathcal{N} (0,\sigma^{2}_{\epsilon^{k}})$, for all $j$. The triple $(\Lambda_{i}(t_{ij}) ,d_{i}^{k},\epsilon^{k}_{ij})$ is a set of independent variables for any choice of $i,j,k$.

A linear model for the intermediate variables $\theta_{ij}^k = \Lambda_{i}({t_{ij}}) +x^k_{i}(t_{ij})^T\beta^k + z^k_{i}(t_{ij})^T d^k_{i}+ \epsilon^{k}_{ij}$ is a particular case of model (\ref{intermediate}).

    \subsubsection{ Link between $\theta_{ij}^k$ and
    the data: the tests functions}

 For subject i, we denote $Y_{ij}^{k}$ the
random variable representing the observation of the $k^{th}$ test on the occasion of
the $j^{th}$ visit at time $t_{ij}$. We will consider the
 cases of ordinal (including binary)
longitudinal data. We consider a test k for which $M_{k}$ ordered values are possible ($m \in
[0,M_{k}-1] $). Observation of
$Y_{ij}^{k}=m$ provides the information that $\theta_{ij}^k$ lies between
two thresholds, that is,
$ Y_{ij}^{k} =   m$ if and only if  $c^{k}_{m} \le  \theta_{ij}^k < c^{k}_{m+1} $, with $c_0=-\infty$ and $c_{M_k}=+\infty$. The test function (which is the function of $\theta_{ij}^k$ that equals $Y_{ij}^{k}$) is in this case a step function. The cut-off points $c^{k}_{m}$ are not known and  must be parameterized
or estimated directly according to the number of possible values
$M_{k}$. Generally we shall represent $c^{k}_{m}$ as a function of parameters $\eta^k$, the dimension of which may be less than $M_{k}-1$ in order to obtain a more parsimonious model: $c^{k}_{m}=\tau^{k}(m,\eta^k), \forall m
\in [1,M_{k}-1]$, where $\tau^{k}(.;\eta^k)$ is a monotone function.

Binary data are simply a special case of ordinal data for which we
only need one cut-off point, $\eta^k_{0}$ for instance. For a binary
test, $  Y_{ij}^{k}  =1_{\{\theta_{ij}^k \ge \eta^k_{0}\}} .$

         \section{Likelihood Inference}

For establishing the likelihood we will first study the distribution of the intermediate variables. Then we establish the likelihood for the case where the tests are ordinal variables as in our application.

       \subsection{Joint distribution of the intermediate variables}

We shall study the distribution of the $Kn_i$ vector $\Theta_i=(\theta_{ij}^k; k=1\ldots, K; j=1,\ldots, n_i)$.
It is to be noted that in equations (\ref{latent}) and (\ref{intermediate})
linearity in the random effects is assumed: this requirement
is important to remain in a Gaussian framework; that is to say  $\Theta_{i}\sim \mathcal{N} (\mu_{i}, \Sigma_{i})$.
 Thus computing the distribution of $\Theta_{i}$ comes down to computing its mean vector $\mu_{i}$ and variance covariance matrix $\Sigma_{i}$.
The expectation can easily be computed since we have:
$${\rm  E }(\theta_{ij}^k)=f( \beta ,x_{i}(t_{ij}) )  + g^k( \beta^k,x^k_{i}(t_{ij}))$$

The variance of $\Theta_i$ is the sum of the variance coming from the latent process $\Sigma_{i,\Lambda }$, the variance of the test specific random effects $\Sigma_{id}$ and the variance of the noise term $\Sigma_{i\varepsilon}$:
$$
\Sigma_{i} = \Sigma_{i\Lambda }+\Sigma_{id}+\Sigma_{i\varepsilon}
 =
\begin{pmatrix}
 \Sigma^{0}_{i\Lambda} & {\dots} & \Sigma^{0}_{i\Lambda} \\
{\vdots} & {\ddots}& {\vdots}\\
 \Sigma^{0}_{i\Lambda} & {\dots} & \Sigma^{0}_{i\Lambda} \\
\end{pmatrix}
+
\begin{pmatrix}
 \Sigma_{id^1} & 0 & 0 \\
0 & {\ddots}& 0 \\
0 & 0 & \Sigma_{id^K} \\
\end{pmatrix}
+\begin{pmatrix}
 \sigma^2_{\varepsilon^1}I_{n_i} & 0 & 0 \\
0 & {\ddots}& 0 \\
0 & 0 & \sigma^2_{\varepsilon^K}I_{n_i} \\
\end{pmatrix},
$$\\
where $\Sigma^{0}_{i\Lambda} = \mathbf{F_{i}}^{T} {\text{~}}  A  {\text{~}} \mathbf{F_{i}} + \Gamma_{i}$,
 and  $\Gamma_{i}$ is the covariance matrix associated with the
Brownian motion:
 $$ \Gamma_{i}= \begin{pmatrix}
t_{i1} & t_{i1} & \dots & t_{i1} \\
t_{i1} & t_{i2} & \dots & t_{i2} \\
\vdots  & \vdots  & \ddots &  \\
t_{i1} & t_{i2} &  & t_{in_{i}} \\
\end{pmatrix},$$
and $\mathbf{F_{i}} = \bigl ( F(\beta,z_{i}(t_{i1})) ,\dots
,F(\beta,z_{i}(t_{in_{i}}))\bigr ) $, a $q \times n_{i}$-matrix, and where
$ \Sigma^k_{id} = \mathbf{G^{k}_{i}}^{T}
{\text{~}}  D^{k} {\text{~}}  \mathbf{G^{k}_{i}}$,
with
$ \mathbf{G^{k}_{i}} = \bigl ( G^{k}(\gamma^k,z^k_{i}(t_{i1})) ,\dots
,G^{k}(\gamma^k,z^k_{i}(t_{in_{i}}))\bigr ) $, a $r_{k} \times n_{i}$ matrix.

\subsection{Identifiability}
Clearly there must be some constraints on the parameters to ensure identifiablity. A thorough analysis is out of the scope of this paper but we give some insight into it. 
We can distinguish three sets of parameters: $\bbeta=(\beta, \beta^k, k=1,\ldots,K)$, $\bgamma=(\gamma, A, \gamma^k, D^k, \sigma_k^2, k=1,\ldots,K)$ and $\boeta=(\eta^k, k=1\ldots,,K)$ and the whole set of parameters is $\balpha=(\bbeta,\bgamma, \boeta)$.
We consider the case of the linear model for the sake of simplicity; in the linear model there is no parameter $\gamma$ nor $\gamma^k$. Clearly in order that $\bbeta$ and $\bgamma$ be identifiable from observation of the $Y_{ij}^k$ they should be identifiable from the observation of $\theta_{ij}^k$.

Let us now look at sufficient conditions for this. In the linear model there is a matrix $\bA$ such that ${\rm E}(\Theta)=\bA \bbeta$. A
necessary and sufficient condition for identifiability of $\bbeta$ is $r(\bA)= dim (\bbeta)$, where $r(\bA)$ is
the rank of $\bA$: this happens if and only if the columns of $\bA$ are linearly independent. A necessary
condition for that is $K\sum n_i \ge dim (\bbeta)$. A sufficient condition of identifiability of $\bbeta$ is:

{\bf C1}: (i) there is no collinearity of the explanatory variables ; (ii) there are no explanatory variable for
one of the equations of the intermediate variable.

Point (i) is common in all linear models.
That C1 is sufficient for identifiability of $\bbeta$  can be seen from the structure of the $\bA$ matrix.

Similarly for the identifiability of $\bgamma$ we consider the condition:

{\bf C2}: (i) There is no random effect for one of the equations of the intermediate variable; (ii) we do not have that all the matrices $F_iF_i^T$ are equal.

For instance if there is no random effect for test $k$ we have:
${\rm var}_{\bgamma}(\theta^k_{i})=F_i^TAF_i+\Gamma_i+\sigma_{\varepsilon^k}I_{n_i}.$
If there was non-identifiability there would exist $\bgamma'\ne \bgamma$ such that ${\rm var}_{\bgamma'}(\theta^k_{ij})={\rm var}_{\bgamma}(\theta^k_{ij})$, which would entail:
$F_i^T(A'-A)F_i=(\sigma_{\varepsilon^k}'-\sigma_{\varepsilon^k})I_{n_i}.$
However the rank of the left-hand side is $q$ while the rank of the right-hand side matrix is $n_i$. So unless $n_i=q$ for all $i$, this equality holds only if $A'=A$ and $\sigma_{\varepsilon^k}'=\sigma_{\varepsilon^k}$. If $n_i=q$ for all $i$, we could solve the equation to find $(A'-A)$ as a function of $F_i$ leading to the additional requirement that $F_iF_i^T$ be the same for all $i$.

As for the identifiability of the whole set of parameters from the observation of  the $Y_{ij}^k$ it is difficult to prove a sufficient condition. There is at least an obvious non-identifiability case that can be detected, and thus avoided. 
For fixed $\bgamma$ the distribution of the $Y_{ij}^k$ depends only on the $c^k_{l}-{\rm E}_{\bbeta}(\theta^k_{ij})$ for $l=1,\ldots, M_{k-1}$, $k=1,\ldots,K$. 
If the model for the cut-off points allows to find $\eta'^k$ such that: $c^k_{l}(\eta'^k)=c^k_{l}(\eta^k)+\Delta$ for $l=1,\ldots, M_{k-1}$, $k=1,\ldots,K$ and if there is an intercept ($\beta_1$) in the equation of the latent process, then the distribution of the $Y_{ij}^k$ specified by $\balpha'$, where $\balpha'$ is defined by $\boeta'$, $\beta_1'=\beta_1+\Delta$ and the other parameters equal to those of $\balpha$, is the same as that specified by $\balpha$. To avoid this non-identifiability case we may for instance give a fixed value to one cut-off value  or the intercept $\beta_1$, a condition we call ``C3''. 

In practice we recommend that conditions C1, C2 and C3 be applied, or analogous conditions since these are
particular cases of constraints that may be put on the three levels of the model.

\subsection{Likelihood}

We will first establish the individual
contribution to the likelihood $\LL_{i}(\balpha)$.
for any subject $i$.
We denote by $y_{ij}^{k}$ the (realized)
observation relative to the $k^{th}$ test on the occasion of the
$j^{th}$ visit at time $t_{ij}$, a realization of $Y_{ij}^{k}$. $\LL_{i}$ is the probability
according to the model of the observed trajectory, that is:
$$ \LL_{i}(\balpha)= P[ Y^{1}_{i1}= y^{1}_{i1},\dots,
Y^{1}_{in_{i}}=y^{1}_{in_{i}}, \dots , Y^{K}_{i1}=
y^{K}_{i1},\dots, Y^{K}_{in_{i}}= y^{K}_{in_{i}} ]
$$

We will now define the sets over which integration will be required.
Let $C_{ij}^{k}$ be the interval relative to observation
$y_{ij}^{k}$ and intermediate variable $\theta_{ij}^k $.
$$
 C^{k}_{ij} = [  c^{k}_{y_{ij}^{k}} , c^{k}_{y_{ij}^{k} + 1} ]
$$

If we define $C_{i}$ the orthant concerning subject $i$,  $ C_{i} = \mathop \otimes \limits^{n_{i},K} _{j=1,k=1}
C^{k}_{ij}$, we obtain for the entire path concerning subject $i$
$$
\LL_{i}(\balpha) =
 P[Y_{ij}^{k} = y_{ij}^{k} ,j=1,\ldots,n_{i}; k=1,\ldots,K  ]
 = P [ \Theta_{i} \in  C_{i} ]
$$

As $ \Theta_{i} \sim  {\cal N}( \mu_{i}, \Sigma_{i}) $, we just need to
integrate the multivariate
normal probability density function $\phi_{(\mu_{i}, \Sigma_{i})}$ over the $C_{i}$ sets:

$$ \LL_{i}(\balpha)= \idotsint \limits_{C_{i}} \phi_{(\mu_{i}, \Sigma_{i})} (\textbf{u})
d\textbf{u}.
$$

Missing values cause no problem because if value at test $k$ at time $t_{ij}$ is missing, the integration set
$C_{ij}^{k}$ for this observation becomes $]-\infty,+\infty[$, so this simply decreases the multiplicity of the integral by one. It is possible to include a truncation condition by
 writing a conditional likelihood. See the application section (4.3) for
 an illustration.
Independence over subjects allows to obtain the likelihood of the sample as  $\LL(\balpha)= \prod _{i=1}^n \LL_i(\balpha)$.

    \subsection{Maximisation algorithm}

The likelihood is difficult to compute since each $\LL_i$ involves a multiple integral, which has to be computed
numerically (see Evans and Swartz, 2000, for a review). However, an advantage of our model is that the integrals that we have to compute are integrals of normal multivariate densities. Efficient techniques exist for this task: in particular  the algorithms proposed by Genz (1992) allow us to compute such integrals up to a
multiplicity of 20. The multiplicity of the integral for computing $\LL_i$ is $K n_{i}$. For instance in
our application we have $K=2$ and $n_i=4$, which leads to a multiplicity of $8$, a feasible problem with the Genz algorithm.

Maximum likelihood estimators can be obtained by using quasi-Newton algorithms. We have considered a Marquardt algorithm (Marquardt, 1963) and an algorithm used by Heddeker and Gibbons (1994) and  Todem, Kim and Lesaffre (2007), in which the Hessian of the log-likelihood is replaced by the estimated variance matrix of the score. This algorithm has been further studied and called ``Robust-variance scoring'' (RVS) algorithm by Commenges et al. (2006). An advantage of the RVS algorithm  is that it needs only first derivatives of the log-likelihood, and the standard errors are obtained from the estimated variance matrix of the score at the maximum. Our experience shows that the RVS algorithm is more than twice as fast as the Marquardt algorithm in our problem.

         \section{Application}
    \subsection{The PAQUID study and the studied sample}
  The proposed approach was applied to the joint modeling of diagnosis of dementia
and a psychometric test, the Mini Mental State Examination (MMSE) (Folstein et al.
1975), using the data of the PAQUID cohort.

 The PAQUID program on cerebral aging is
based on a large cohort randomly selected in a population of subjects
aged 65 years  or older, living at home in two administrative areas of southwest
France (Gironde and Dordogne). Our analysis bears on the first eight years of the follow-up of this study. In addition to the initial visit, subjects were seen approximately after one, three, five and eight years in Gironde and after three, five and eight years in Dordogne; the successive visits are denoted by T0, T1, T3, T5 and T8. At each visit the MMSE was measured and
diagnosis of dementia was made by neurologists, based on the DSMIII-R criteria (for details see Letenneur et al., 1999). We will use the first five years to fit the model and the eight-year follow-up to assess the predictive ability of our model.

Our sample was composed only of women who were not demented at the initial visit. It is safer to analyse men and women separately because the dynamics of ageing seems to be quite different between the genders (see Commenges et al., 2004). Because there are more women than men in the PAQUID sample we chose to focus on women. We introduced the condition of being non-demented at the initial visit because it is doubtful that the PAQUID sample is representative of the whole population (demented and non-demented): demented subjects are often institutionalized. The condition of being non-demented at entrance must be taken into account in the likelihood (see section 4.3).
 At the initial visit there were five cases which, although not diagnosed as demented, obtained a MMSE score of zero (this can be seen on Figure 1): these subjects had cognitive impairment due to other causes than dementia (stroke, psychiatric illness); we have chosen to keep them in the sample.

We thought that the evolution of cognitive ability may be strongly affected by dementia and it was not our aim to describe this evolution; in consequence, further observations of the MMSE after diagnosis of dementia were not taken into account. This artificial right-censoring is ignorable: the reason is that it is done on the basis of an observed variable included in the model and this can be proved using results of Commenges and G\'egout-Petit (2005).

Finally, our study sample was composed of 2131 women aged 65 years or older and who were not demented at the initial visit.
During the 5-year follow-up we had 5622 observations of the MMSE. We had also 5742 assessments of the demented
status; among them, 126 were diagnoses of dementia.

    \subsection{The model applied to the PAQUID sample}
\subsubsection{The explanatory variables}

The different components of the model we developed may depend on educational level and a variable indicating whether the test was administered for the first time (to take into account a possible practice effect): educational level has been shown to be a risk factor of dementia (Letenneur et al., 1999) and a practice effect of the MMSE has been found (Jacqmin-Gadda et al., 1997). Moreover, there has been debate about the necessity of correcting the MMSE for educational level in order to determine cognitive impairment, a prognostic factor of dementia.

The most difficult problem is to define what time is in our model. Since we wish to relate cognitive decline to age it is natural to determine a time-scale for each subject closely related to age. We could consider that the time that is relevant for a subject is the time elapsed since her birth, that is, age. However, in this model we do not wish to model the evolution of cognitive ability from birth (we would have to develop a much more complicated model) but only the decline of cognitive ability from an age at which we think that this phenomenon may start for a non-negligible fraction of the population. We took as origin the age of 65 for the two following reasons: (i) we have observations from age 65, making it awkward to take a later origin, which would lead to negative times: particularly in a non-stationary (due to the Brownian motion) and non-time-homogeneous (due to the non-linearity in $t$) model this would not make sense; (ii) we have tried earlier time origins but this yielded lower likelihoods.

  Educational level is represented by the binary variable that we will denote by ${\Ed}_i$ so that
$ \Ed_i =  1$ if subject i has obtained a primary school diploma and $0$
if not.  Practice effect, denoted by ${\Fp}_i$, is defined as:
    ${\Fp}_i (t) =1$ for $t \leq t_{i1}$ and
 ${\Fp}_i (t) =0$ for $t > t_{i1}$.

For clarity of interpretation we will describe the model directly in terms of $t$, $\Ed_i$ and $\Fp_i(t)$ rather than using the general notations.

    \subsubsection{The latent process}
 In this application of our model, the latent process represents
cognitive ability: diagnosis of dementia and MMSE will be considered as indirect measurements of it. The latent process is defined by equation (1)  in which we specify $f(.,.)$ as:
$$f(\beta , x_{i}(t)) = (\beta_{1} + \beta_{2}\Ed_i ) +
(\beta_{3} + \beta_{4}\Ed_i) t^{\beta_{5}},$$
 As for the function $F(.,.)$ we tried:

$ F(\gamma,z_{i}(t)) a_{i} =
\begin{pmatrix}
 1 , t^{\gamma_{1}}
\end{pmatrix}
\begin{pmatrix}
\ a_{1,i} \\
\ a_{2,i} \\
\end{pmatrix}
=  a_{1,i}  + a_{2,i}  t^{\gamma_{1}}.$
It was natural to assume $\gamma_1=\beta_5$,
that is, there is a vector of random effects $\mathbf{a_{i}}$
of size $q=2$ bearing on the intercept $\beta_{1}$ and the slope
$\beta_{3}$. However the algorithm failed to converge when we tried to estimate the two variance parameters and the correlation coefficient of the two random effects, probably due to the presence of the Brownian motion.
The algorithm converged if we assumed a diagonal variance matrix for  $a_{i}$:
$ A=\begin{pmatrix}
 \sigma_{a_{1}}^{2} & 0 \\
 0 & \sigma_{a_{2}}^{2} \\
\end{pmatrix}$. We also tried a simpler model with only one random effect obtained with the $F(\gamma,z_{i}(t)) a_{i} =a_i$; since this simpler model gave nearly the same result, we present this simpler model in the following.
For this model the latent process is defined as:
\begin{equation} \Lambda_i(t)=(\beta_{1} + \beta_{2}\Ed_i ) +
(\beta_{3} +  \beta_{4}\Ed_i)t^{\beta_{5}}+a_{1,i} +W^i(t).\end{equation}

    \subsubsection{Observation equations.}
In this application, we jointly model the diagnosis of dementia
and the MMSE score, so that $K=2$: the first ``test'' (k=1) is diagnosis of dementia and this is a binary variable; the second ``test'' (k=2) is the MMSE which has 31 values. The specification of the equations for the intermediary variables is guided by interpretability and identifiability issues.

We have introduced a random effect in the model of the intermediate variable $\theta^1_{ij}$ for  diagnosis $(k=1)$.  In formula (\ref{intermediate}) we took $g_i^{1}(\beta^1,x^1_{i} (t_{ij}))=0$ and
$G_i^{1}(\gamma^1,z^1_{i} (t_{ij}))=1$; there was one random effect $d_{i}^{1} \sim N(0,\sigma^{2}_{d^1} )$. This random effect  makes it possible that subjects with a low latent process are not diagnosed demented; this may happen because some subjects have always had low cognitive ability not linked to a neurodegenerative process. We did not introduce additional error term, that is to say $\sigma^2_{\varepsilon ^1}=0$, nor explanatory variables (thus satisfying condition C1 in section 3.2). Thus the intermediate variable for dementia is:
\begin{equation}\theta^1_{ij}=\Lambda_i(t_{ij})+d^1_i.\end{equation}
For relating this variable to the diagnosis of dementia (which means defining the ``test function'') we just need one cut-off value given by the parameter $\eta_0$:
$Y^1_{ij}=1 $ if and only if $\theta^1_{ij}\le \eta_0$. Our notation here for the parameters $\eta$ differs slightly from the general case: we use $\eta_0$ for dementia and $\eta_{1}$, $\eta_{2}$ and $\eta_{3}$ for the MMSE, the meaning of which is explained below.

As for the MMSE $(k=2)$ we took
into account both the practice effect and the specific impact
of educational level on MMSE. The practice effect is only
located on the first visit $(j=1)$ and we introduced an interaction with
educational level (meaning that the practice effect may not be the same for subjects with or without a primary school diploma). Thus in formula (\ref{intermediate}) we took
$g_{i}^{2}(\beta,x_{i} (t_{ij}))=
\beta^2_{1}\Ed_i +  \beta^2_{2} \Fp_i(t_{ij}) + \beta^2_{3}\Ed_i\times \Fp_i(t_{ij}).$ No specific random effect  was
introduced in the MMSE equation (condition C2), so $ G_i^{2}(\gamma^2,x_{i} (t_{ij}))= 0 $. There was, however, an error term of variance $\sigma^2_{\varepsilon ^2}$. Thus, the intermediate variable for MMSE was:
\begin{equation} \label{MMSE} \theta^2_{ij}=\Lambda_i(t_{ij})+\beta^2_{1}\Ed_i +  \beta^2_{2} \Fp_i(t_{ij}) + \beta^2_{3}\Ed_i\times \Fp_i(t_{ij})+\varepsilon ^2_{ij}.\end{equation}

MMSE takes values between $0$ and $30$, so we have $M_{2}=31$.  It is judicious to use a model for the family of cut-off points $c^{2}_{m}=\tau^{2}(m,\eta)$ which is more parsimonious than considering all the cut-off values as parameters. We have $c^{2}_{M_{2}}=+\infty$ and $c^{2}_{0} = - \infty$ and for satisfying condition C3 we fixed $c^{2}_{M_{2}-1}$ arbitrarily at the value $c^{2}_{M_{2}-1}=40$. There is no reason that the MMSE scale be linear with respect to the latent process scale so we used the following model yielding unequally spaced cut-off points:
$c^{2}_{m}=  40 - \eta_{1} (M_{2}-1-m)^{\eta_{2}}$. We limited this power model to $m \in [1,M_{2}-3]$ and we gave an independent parameter $\eta_{3}$ for $c^{2}_{M_{2}-2}$, which made it possible to improve the fit as compared to extending the above model up to $M_{2}-2$. Thus our model for the test function for MMSE involves three parameters: $\eta_{1}$, $\eta_{2}$ and $\eta_{3}$.

    \subsection{The likelihood for the application}

We computed the likelihood according to section 2.
We also had to include the selection condition mentioned in section 4.1:
 since only non-demented subjects were included, the likelihood is conditional on $\{\theta^1_{i1} > \eta_{0}\}$ (the event that subject $i$ is not diagnosed demented at
initial visit ${t_{i1}}$); the conditional likelihood for subject $i$ is $\LL_i/ P(\theta^1_{i1} > \eta_{0})$.
We obtain from the model: $ \theta^1_{i1} \sim
\mathcal{N}
\begin{pmatrix}
  & f(\beta , x_{i}(t_{i1}))&,&  \Sigma_{i}( 1, 1 ) & \\
\end{pmatrix}
 $, so that we have:
$$ P( \theta^1_{i1} > \eta_{0} ) =
\Phi\Bigl ( \frac {f(\beta , x_{i}(t_{i1})) - \eta_{0} }
{\sqrt{\sigma_{a_{1}}^{2}
+ t_{i1} + \sigma_{d^1}^{2} }} \Bigr ).$$

The likelihood was maximized using the RVS algorithm described in section 3.3.

    \subsection{A Simulation}
In order to demonstrate the ability of our algorithm to maximize such a complex likelihood we tried it on a simulated data set.
    We generated a sample of size $n=2131$ with the same age
    distribution at the initial visit and the same proportion of
    educated and non-educated subjects as in the real data sample from the PAQUID study. We generated  4 visits as in the real data set, the initial visit and visits after one, three and five years.
    The values of the 16 parameters were taken equal to the values estimated in the real data set.  We took as starting values: $\beta_2=\beta_3=\beta_4=\beta^2_1=\beta^2_2=\beta^2_3=0$; $\beta_1=38.5$;  $\beta_5=1$; $\eta_0=30$; $\eta_1=\eta_2=1$; $\eta_3=39$; $\sigma_{a_1}=10^{-5}$; $\sigma_{d^1}=\sigma_{\varepsilon^2}=10$. The algorithm converged in $19$ iterations. The results are given in Table 1. We see that the estimated values are reasonably close to the target values and that the .95 confidence intervals include these values. The algorithm converged toward the same point from different starting values. We also verified the quality of the inverse hessian for giving estimates of the variances of the estimators of the parameters by checking a reasonable agreement between some Wald tests and likelihood ratio tests. On the whole, the algorithm seems to be reliable.

    \subsection{Model estimated from the PAQUID data}

 The values of the parameters estimated from the PAQUID sample are shown in Table 2.
 As expected there is a significant mean trend of decrease of global cognitive ability (see $\beta_3$) with a shape not far from a quadratic form (see  $\beta_5$). There is a significant heterogeneity around the intercept (see $\sigma_{a_1}$). The significant random effect for dementia ($\sigma_{d^1}$) means that some subjects are not diagnosed demented at repeated visits in spite of low cognitive ability.

The value of 0.58 for parameter $\eta_{2}$ indicates that a difference of one point of MMSE corresponds to a
larger difference in cognitive ability for high cognitive level than for a low one; in other words, the sensitivity
of MMSE is better for low level than for high level; this is graphically illustrated in Figure 2 which displays
a grid of the cut-off values making it visible that a larger difference in latent process (or rather intermediate
variable) values is necessary to make one point of difference for the MMSE for higher rather than for lower level.
This is reminiscent of the mixed linear model applied by Jacqmin-Gadda et al. (1997) to the square-root of 30
minus MMSE (in fact the number of errors).

 In order to assess the degree of realism of our model for the MMSE we
computed the expected numbers of subjects having score $m$ at the MMSE at T0: this was achieved by computing for
each subject the probability of having score $m$ and summing the 2131 probabilities. The computation of these
probabilities was carried out with the estimated model, taking into account the ages, educational levels and the practice effect, as well as the different variability terms and the use of formulae similar to that used for the
prediction in section 4.6.  Figure 1 compares the histograms of observed MMSE scores with the histogram of
expected numbers; it can be seen that they are quite similar. There is a slight discrepancy at scores 22 and 21:
this artefact is due to the screening design for diagnosing dementia in the PAQUID study at T0 which used the
threshold 24 and which probably led interviewers to put 22 or 21 rather than 24 for some subjects (to trigger the visit of a neurologist).

We can make an approximate link between the threshold for dementia $\eta_0$ and values at the MMSE. Taking zero values for the random effect for dementia and errors for the MMSE, the value of the threshold approximately corresponds to MMSE= 19 and MMSE= 21 for low  and high educational levels respectively. (The value 19 is found as follows. For a subject with low educational level we have from (6): $\theta^2_{ij}=\Lambda_i(t_{ij})$ and ${\rm E}(\theta^1_{ij})=\Lambda_i(t_{ij})$; thus if we consider subjects for which ${\rm E}(\theta^1_{ij})=\eta_0$ they have $\theta^2_{ij}=\eta_0$; the corresponding value $m_0$ of the MMSE score satisfies the equation $\eta_0=40-\eta_1(30-m_0)^{\eta_2}$).

 Our model allows us to distinguish the effect of educational level on the latent cognitive ability on the one hand and on the MMSE score on the other. Educational level has a significant effect ($\beta_2$) on the intercept of the cognitive ability process, but not on the slope ($\beta_4$); there is a highly significant effect of educational level ($\beta^2_1$) for the MMSE. To sum up, (because of the positive $\beta_2^2$) subjects with high educational level tend to have higher MMSE than subjects with low educational level, for the same value of the latent process (true cognitive ability), leading to a diagnosis of the former as demented at higher MMSE levels than the latter; on the other hand (because of the positive $\beta_2$) subjects with high educational level tend to have higher value of the latent process than subjects with low educational level, leading to a lower rate of diagnosis of dementia for the former as compared to the latter. Finally, there is a significant effect of practice effect ($\beta^2_2$)(subjects have a lower MMSE at the first visit than what would be expected); the interaction of practice with educational level ($\beta^2_3$) is not significant.

Several features of these results can be best illustrated by a graphic.
Figure 2 displays, in the latent process scale, both the grid of the cut-off values for the MMSE (horizontal dotted lines) and the threshold for diagnosis of dementia (the horizontal crosses line at $\eta_0=24.41$). It also displays the expected value of the latent process of cognitive ability for subjects with low and high educational level (the curve for low educational level starts at the value of the intercept $\beta_1=32.90$). The curves are approximately parallel and the curve for low educational level below; this explains that a larger incidence of dementia has been observed in this group (Letenneur et al., 1999). We can see that the decline of this expected value is very slow near the age of 65 and accelerates for older ages for both low and high educational levels. This is rather in agreement with normative values which have been established in the United States (Crum et al., 1993) and in France (Lechevallier-Michel et al., 2004) although the results can not be compared directly: one main difference is that normative values exclude demented subjects; another difference is that we model the practice effect. Figure 2 also shows the dispersion for the values of the latent process by showing a region in which $95\%$ of the values for low educated subjects lie at each age. The lowest bound curve (dashed line) crosses the threshold value (around 75) and so, it is graphically apparent that a growing number will be diagnosed demented with older age.

Moreover Figure 2 illustrates the effect of educational level on values of the MMSE (for a given value of the latent process), as well as the practice effect on MMSE scores. It displays the expected values of intermediate variables for MMSE ($\theta^2_{ij}$) for subjects with low and high educational levels entering at 75 in the study and seen one, three, five and eight years after. In our model these expected values are equal to the expected value of the latent process for subjects with low educational level (the stars) except for the first visit where the value is lower due to the practice effect: this is because if $\Ed_i=0$ and $\Fp_i=0$ we have from formula (\ref{MMSE}) $\theta^2_{ij}=\Lambda_i(t_{ij})+\varepsilon ^2_{ij}$, so that ${\rm E}(\theta^2_{ij})={\rm E}[\Lambda_i(t_{ij})]$. As already mentioned, there is a grid indicating the values of the MMSE obtained as a function of the intermediate variable. For instance a subject with low educational level who has her intermediate variables equal to the expectations and entering at 75 at T0 would have MMSE values 24, 25, 25, 24 and 23 at T0, T1, T3, T5 and T8 respectively. The expectations of the intermediate variables for subjects with high educational level are higher than the expected value of the latent process for the same time. The results illustrated in this figure, contribute to the debate regarding the possible correction of the MMSE to take the educational level into account and regarding the effect of educational level on dementia. It appears that educational level has an effect on global cognitive ability (our latent process), and thus on dementia, but also has a specific effect on MMSE.

    \subsection{Prediction of dementia diagnosis}

    The model may be used for predicting diagnosis of dementia for
    subject $i$ at time $t_{i,n_i+1}$, given the MMSE values at the successive
    visits $(1,\ldots,n_{i})$ and given that subject $i$ has not
    been diagnosed demented up to visit $n_{i}$.
The information that we have up to visit $n_{i}$ is summarized by the event $\Theta_i \in C_i$.
The probability that subject $i$ is diagnosed demented at $t_{i,n_i+1}$ is
$$p_{i}=P[\theta^1_{i,n_i+1} \le \eta_{0}|\Theta_{i} \in C_{i}]= \frac
{P[(\theta^1_{i,n_i+1} \le \eta_{0}) \cap (\Theta_{i} \in C_{i})]}
{P[\Theta_{i} \in C_{i}]}.$$
This expression is not affected by the condition of not being
diagnosed demented up to visit $n_{i}$ as the corrective
conditional probability cancels out in the ratio.
In order to compute the numerator we need the joint distribution of $\theta^1_{i,n_i+1}$ and $\Theta_i$. This is a normal distribution with
expectation:
$$
\mu_{i}^{*}=\begin{pmatrix}
  \mu_{i} \\
   E[\theta^1_{i,n_i+1}]\\
\end{pmatrix}
=\begin{pmatrix}
  \mu_{i} \\
  f( \beta ,x_{i}(t_{i,n_i+1}) ) \\
\end{pmatrix},$$
and variance matrix $\Sigma_i^*$ formed by the block $\Sigma_i$ augmented by the correlation between $\theta^1_{i,n_i+1}$ and $\Theta_i$ and the variance of $\theta^1_{i,n_i+1}$. These are given by:
$${\rm cov}(\theta^1_{i,n_i+1},\theta^1_{ij}) =
 \sigma_{a_{1}}^{2}+ t_{ij} + \sigma_{d}^{2}, \mbox { for } j=1,\ldots,n_i+1;$$
$${\rm cov}(\theta^1_{i,n_i+1} ,\theta^2_{ij}) =
 \sigma_{a_{1}}^{2}
+ t_{ij}, \mbox { for }  j=1,\ldots,n_i.$$

    We selected subjects that had not been diagnosed demented up
    to visit T5 and who had been seen at T8: $N=1187$ subjects satisfied these criteria.
    We computed their individual probabilities $p_i$ of being
    diagnosed demented at visit T8, using the values of the parameters $\theta$ estimated from the follow-up up to five years.
One aim was to predict the number $N_d$ of subjects diagnosed demented at T8: a natural predictor is the expectation of $N_d$ (conditional on information up to T5)
 which  is $\sum_{i=1}^N p_i$. We found $\hat N_d=46.6$. A predictive interval can be computed using the fact that ${\rm var} N_d=\sum_{i=1}^N p_i(1-p_i)$ and treating $N_d$ as approximately normal; we found that the $95 \%$ predictive interval was $[34.1;59.2]$.
We observed 56 new diagnoses at T8, a number inside the predictive interval.

Another way to assess the predictive ability of our model for diagnosis of dementia at T8 was to consider the $p_i$'s as quantitative values on which a classification as positive or negative could be made according to a cut-off value, as in the theory of diagnostic tests. Sensitivity and specificity can be computed for each cut-off value and the ROC curve relates sensitivities and specificities for the different cut-off values. Figure 3 gives the ROC curve for  our prediction of dementia diagnosis. In particular, the area under the ROC curve is a summary measure of performance of the test. The area under the ROC curve of our model is $0.82$, a rather good value.

         \section{Conclusion}
We have developed a general model for multivariate longitudinal ordinal data. It could be easily extended to
include continuous data: we could use for test $k$ a continuous function $h_k(.)$ : $Y_{ij}^{k}  =
h_k(\theta_{ij}^k).$ Such a test function could be chosen in a family of functions depending on a parameter
$\eta^k$. For instance Proust et al. (2006) in an analogous problem have chosen the family of beta cumulative
distribution functions indexed by two parameters.

 When modeling cerebral ageing one would also  have to model death: joint modeling of dementia and death has been achieved by the use of an illness-death model (Joly et al., 2002; Commenges et al., 2004) but cognitive ability was not modeled. It is not possible to rigorously treat the joint occurrence of diagnosis of dementia, psychometric tests and death with existing models. However, approximate inference can be made by considering death as censoring, as has been done in this paper.

Our model is useful for jointly modeling psychometric tests and diagnosis of dementia but could be applied to other epidemiological contexts.

\section*{References}
\noindent
\setlength{\parindent}{-8mm}

{   Aguilar O, Huerta G, Prado R, and West M} (1999). {  Bayesian inference on latent structure in time series (with discussion)},  Bayesian Statistics 6, J.O. Berger, J.M. Bernardo, A.P. Dawid, and A.F.M. Smith (Eds.), Oxford University Press.

{  Chib S and Greenberg E} (1998). Analysis of multivariate probit models. {  Biometrika}, {  85}: 347-361.

{  Commenges D and G\'egout-Petit A} (2005). Likelihood inference for incompletely observed stochastic processes: general ignorability conditions. {  arXiv:math.ST/0507151}, http://arxiv.org/abs/math/0507151.

{  Commenges D, Joly P, Letenneur L and Dartigues JF} (2004). Incidence and
prevalence of Alzheimer's disease or dementia using an
Illness-death model. {  Statistics in Medicine} {   23}: 199-210.

 Commenges D, Jacqmin-Gadda H, Proust C and Guedj J (2006). A Newton-Like algorithm for likelihood maximization: the robust-variance scoring algorithm.{\em arXiv:math.ST/0610402}, http://arxiv.org/abs/math/0610402.

{  Crum RM, Anthony JC, Bassett SS and Folstein MF} (1993). Population-based norms for the mini-mental state examination by age and educational level, {  JAMA}, {  18}: 2386-2391.

{  Dunson DB} (2003) Dynamic latent trait models for multidimensional longitudinal data. {  Journal of the American Statistical Association}, {  98}: 555-563.

{  Evans M and Swartz T} (2000). Approximating integrals via Monte Carlo and deterministic methods. {\it Oxford University Press; Oxford}.

{  Faucett CL and Thomas DC} (1996). Simultaneously modeling censored survival data and repeatedly measured covariates: a Gibbs sampling approach. { Statistics in Medicine} {  15}: 1663-1685.

{  Folstein MF, Folstein SE and McHugh PR} (1975).  Mini-Mental State. A practical
method for grading the cognitive state of patients for the clinician.
{  Journal of Psychiatric Research} {  12}: 189-98.

{  Genz A} (1992). Numerical computation of the multivariate normal probabilities. {  Journal of Computational and Graphical Statistics} {  1}: 141-150.

{  Hashemi R, Jacqmin-Gadda H and Commenges D} (2003). A latent process model for joint modeling of events and marker. {   Lifetime Data Analysis} {  9}: 331-343.

{  Hedeker D and Gibbons R} (1994). A random-effects ordinal regression model for multilevel
   analysis. {  Biometrics} {  50}: 933­944.

{  Henderson R, Diggle  P and Dobson A} (2000). Joint modeling of longitudinal
measurements and event time data. {\it Biostatistics} {  1}: 465-480.

{  Jacqmin-Gadda H, Commenges D and Dartigues JF} (2006). Random changepoint model for joint modeling of cognitive decline and dementia, {   Biometrics}, {  62}: 254-260.

{  Jacqmin-Gadda H, Fabrigoule C, Commenges D and Dartigues JF} (1997).  A 5-year longitudinal study of the mini-mental state examination in normal aging. { American Journal of Epidemiology}  {  145}: 498-506.

{  Joly P, Commenges D, Helmer C and Letenneur L} (2002). A penalized
likelihood approach for an illness-death model with interval-censored
data: application to age-specific incidence of dementia. {Biostatistics} {  3}: 433-443.

{  Lechevallier-Michel N, Fabrigoule C, Lafont S, Letenneur L, Dartigues JF} (2004). Normative data for the MMSE, the Benton visual retention test, the Isaacs's set test, the digit symbol substitution test and the Zazzo's cancellation task in subjects over the age 70: results from the PAQUID Study
{  Revue de Neurologie}  {  160}: 1059-1070.

{  Lee MLT, Degruttola V and  Schoenfeld D} (2000).
   A model for markers and latent health status.      { Journal of Royal Statistical Society: Series B} {  62}: 747-762.

{  Letenneur L,  Gilleron V, Commenges D,Helmer C, Orgogozo JM and
Dartigues JF} (1999).  Are sex and educational level independent predictors
of dementia and Alzheimer's disease? Incidence data from the PAQUID
project. {\it Journal of Neurology Neurosurgery and Psychiatry}
{  6}: 177-183.

{  Marquardt D} (1963). An algorithm for least-squares estimation of nonlinear
parameters. { SIAM Journal of Applied Mathematics} {  11}: 431-441.

{  McCullagh P and Nelder JA}, (1989). Generalized linear models. Chapman \& Hall/CRC, Boca Raton.

{  Proust C, Jacqmin-Gadda, H, Taylor JMG, Ganiayre J and Commenges D} (2006). A nonlinear model with latent process for cognitive evolution using multivariate longitudinal data. Biometrics,  {  62}: 1014-1024.

{  Todem D, Kim KM and Lesaffre E} (2007).
Latent-Variable Models for Longitudinal Data with Bivariate Ordinal Outcomes. {  Statistics in Medicine}, in press.

{  Whitmore GA, Crowder MJ and Lawless JF} (1998). Failure inference from a
marker process based on a bivariate Wiener model. { Lifetime Data Analysis} {  4}: 229-251.

{  Wulfsohn MS and  Tsiatis AA} (1997). A joint model for survival and longitudinal data measured with error. { Biometrics} {  53}: 330-339.


\newpage
\begin{table}
\caption{A simulation mimicking the PAQUID study example.}
\begin{center}
\begin{tabular}{|c|rrr|}
\hline
Parameters  & Targets & Estimates  & St. Dev. \\
\hline
$\beta_{1}$ & 32.90 & 32.51 & 0.36 \\
$\beta_{2}$ & 2.34  & 3.09 & 0.46\\
$\beta_{3}$ & -0.022  & -0.017 & 0.006\\
$\beta_{4}$ & 0.0013 & 0.02 & 0.13\\
$\beta_{5}$ & 1.84  & 1.91 & 0.10\\
$\beta_{1}^{2}$ & 1.69  & 1.41 & 0.35\\
$\beta_{2}^{2}$ & -1.65 & -1.53 & 0.15\\
$\beta_{3}^{2}$ &  0.29  & 0.25 & 0.17\\
$\eta_{0}$ & 24.41 & 24.38  & 0.60\\
$\eta_{1}$ & 3.93  & 3.94 & 0.16 \\
$\eta_{2}$ & 0.58  & 0.58 & 0.01\\
$\eta_{3}$ & 36.64 & 36.52 & 0.15\\
$\sigma_{a_{1}}$ & 2.04  & 2.10 & 0.21\\
$\sigma_{D^{1}}$ &  2.68  & 2.49 & 0.18\\
$\sigma_{\varepsilon^{2}}$ & 2.55 & 2.59 & 0.11\\
\hline

\end{tabular}
\end{center}
\end{table}
\newpage
\begin{table}
\caption{Results from the analysis of the five-year follow-up of the PAQUID study}
\begin{center}
\begin{tabular}{|l|rr|}
\hline
Parameters  & Estimates  & St. Dev. \\
\hline
$\beta_{1}$: intercept for $\Lambda$ & 32.90 & 0.41\\
$\beta_{2}$: effect of education on intercept& 2.34  & 0.55\\
$\beta_{3}$: slope of $\Lambda$ & -0.022 & 0.008 \\
$\beta_{4}$: effect of education on slope & 0.0013 & 0.0018\\
$\beta_{5}$: power of $t$ & 1.84  & 0.11\\
$\beta_{1}^{2}$: effect of education on MMSE& 1.69  & 0.45\\
$\beta_{2}^{2}$: practice effect for MMSE & -1.65 & 0.17\\
$\beta_{3}^{2}$: interaction education x practice effect &  0.29  & 0.20\\
$\eta_{0}$: threshold for dementia & 24.41 &  0.65\\
$\eta_{1}$: multiplicative factor for the cut-off model of MMSE & 3.93  &  0.19\\
$\eta_{2}$: power for the cut-off model of MMSE & 0.58  &  0.006\\
$\eta_{3}$: value of $c_{29}$& 36.64 &  0.17\\
$\sigma_{a_{1}}$: variance of the random effect for intercept& 2.04  & 0.21 \\
$\sigma_{D^{1}}$: variance of the random effect for dementia &  2.68  &  0.20\\
$\sigma_{\varepsilon ^2}$: variance of error in the intermediate equations for MMSE& 2.55   & 0.13\\
\hline

\end{tabular}
\end{center}
\end{table}

\newpage

\begin{figure}
\begin{center}
\includegraphics{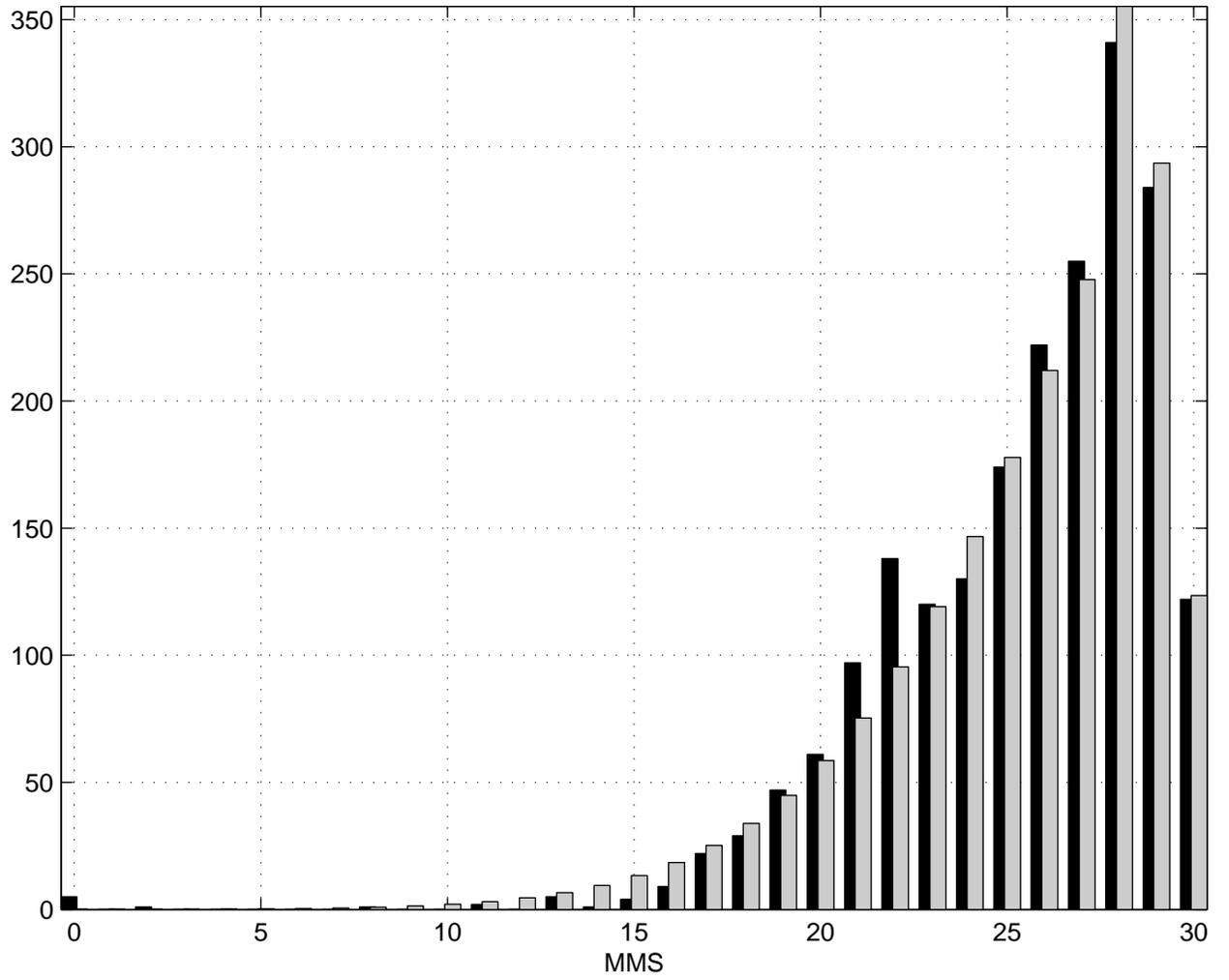}
\caption{Histogram of the MMSE score at the initial visit. Black: observed histogram; grey: expected numbers.}
  \label{}
\end{center}
\end{figure}

\newpage

\begin{figure}
\begin{center}
\includegraphics{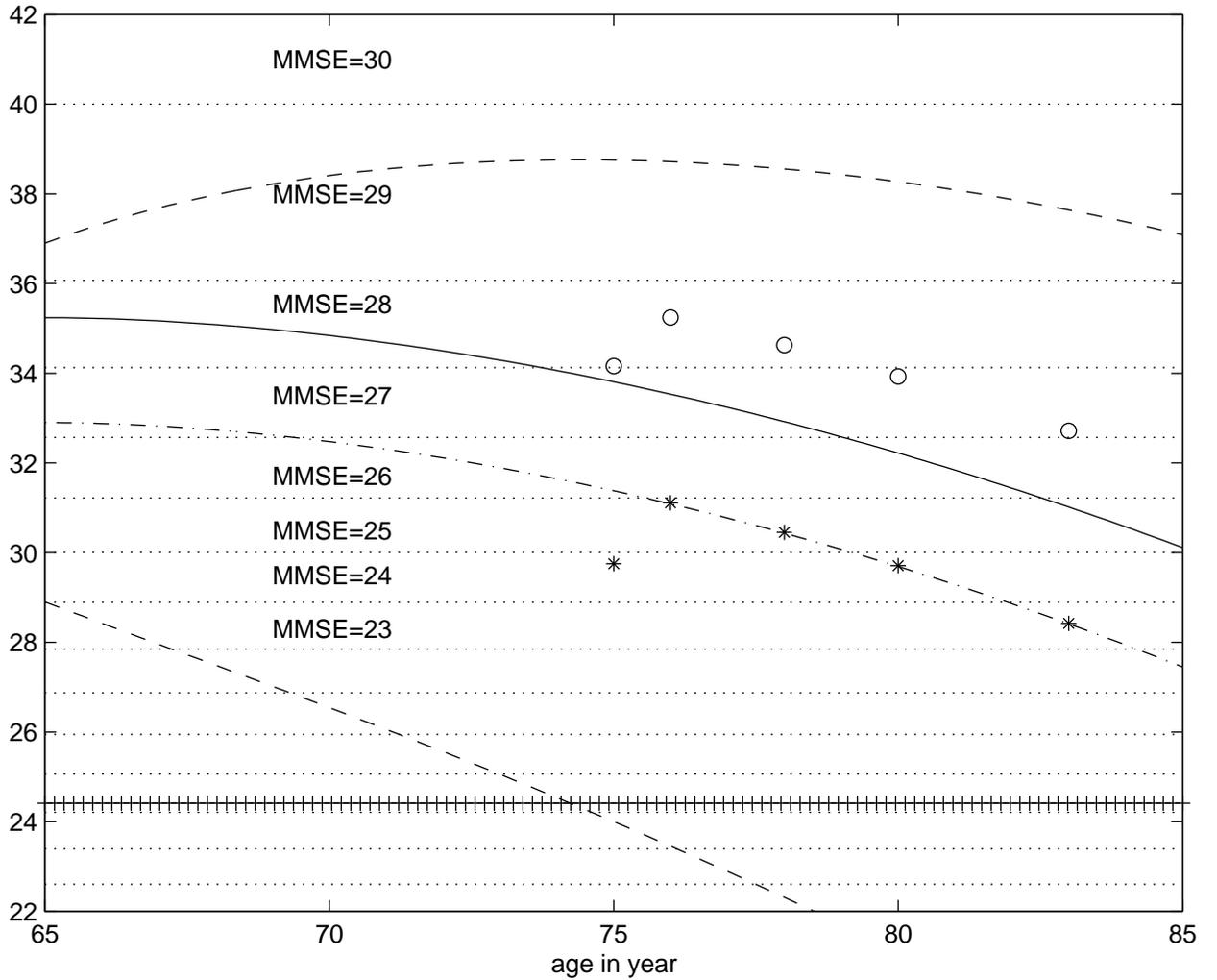}
\caption{Mean evolution of the latent process based on the follow-up of five years in the PAQUID study for low (dashed line) and high (plain line) educational level; the band (delimited by the dashed lines) shows a region where $95\%$ of the values for low educated subjects lie; horizontal line with crosses is the threshold value for dementia; expected intermediate variables for subjects of low (stars) and high (open circles) educational level entering at 75 years in the study and seen at T0, T1, T3, T5 and T8; the grid shows the values of the MMSE obtained for specific values of the intermediate variable. }
  \label{}
\end{center}
\end{figure}

\newpage

\begin{figure}
\begin{center}
\includegraphics{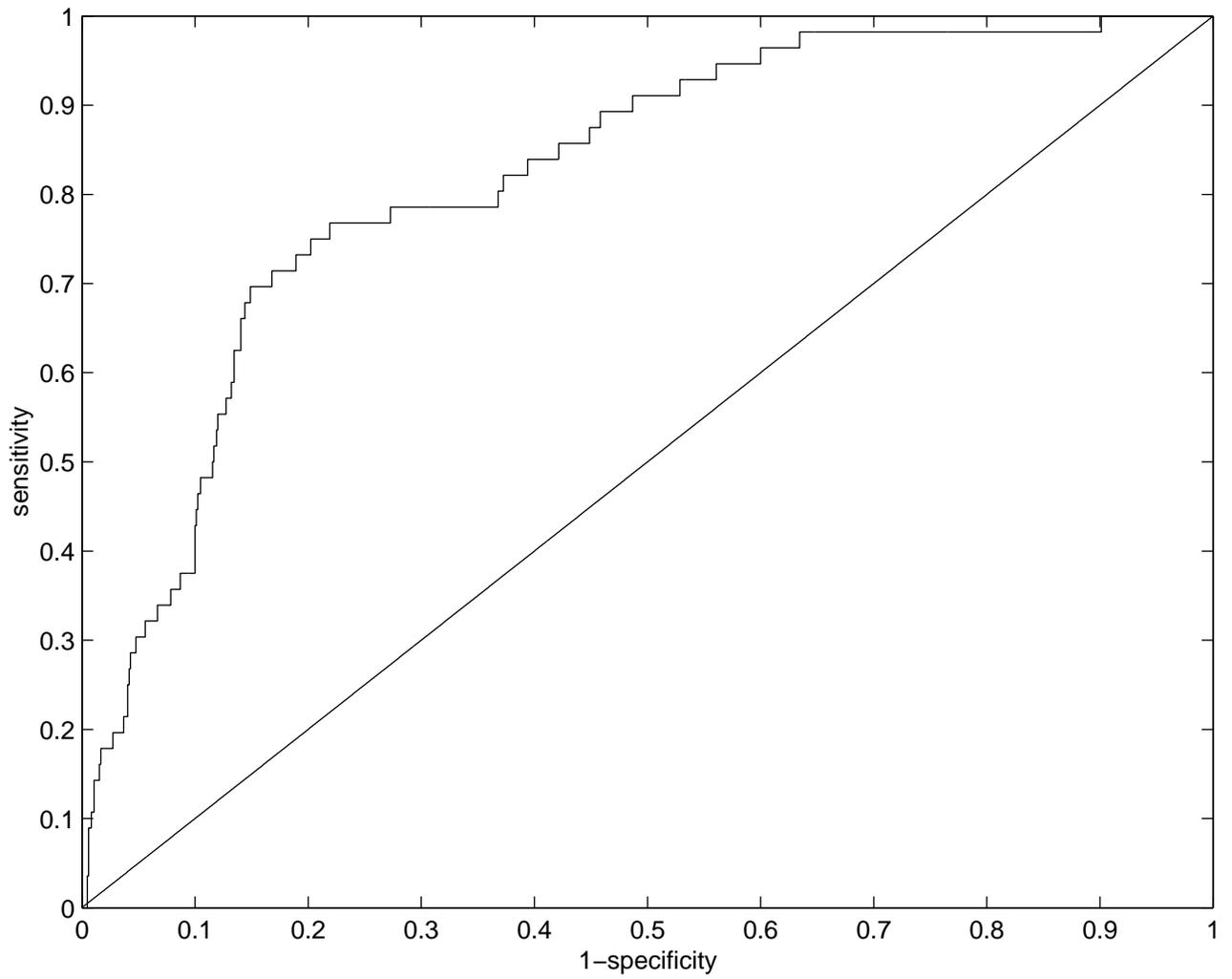}
  \caption{ROC curve showing the ability of the model to predict dementia at the eight-year visit based on the follow-up of five years in the PAQUID study}
  \label{}
\end{center}
\end{figure}

 \end{document}